\documentclass[12pt,reqno]{amsart}
\usepackage{fullpage}
\usepackage{mypack}

\usepackage{tikz}
\usetikzlibrary{matrix,arrows,decorations.pathmorphing}
\usepackage{tikz-cd}
\usetikzlibrary{arrows}

\usepackage{amsmath}
\usepackage{mathrsfs}

\usepackage{epsfig}
\usepackage{amsmath,amsfonts,amsthm,amscd, amssymb}
\usepackage{latexsym}
\usepackage{enumerate}
 
\usepackage{chngcntr}
\makeatletter
\let\c@equation=\c@subsubsection

\makeatother

\title{Spherical posets from commuting elements}
\date{}
\address{Department of Mathematics, University of Western Ontario, London ON  N6A 5B7}
\email{cokay@uwo.ca}
\author{C\.{I}han Okay}

\begin{document}
  \maketitle 
 
\begin{abstract}
In this paper we study the homotopy type of the partially ordered set of left cosets of abelian subgroups in an extraspecial $p$--group. We prove that the universal cover of its nerve  is homotopy equivalent to a wedge of $r$--spheres where $2r\geq 4$ is the rank of its Frattini quotient. This determines the homotopy type of the universal cover of the classifying space  of transitionally commutative bundles as introduced in \cite{AG15}. 
\end{abstract}  
\section{Introduction}
The motivation behind this work is initiated by our interest in a certain filtration of the classifying space $BG$ of a group $G$ introduced in \cite{ACT12}. This filtration consists of subspaces  $B(q,G)\subset BG$ for $q\geq 2$ defined as geometric realizations of certain simplicial spaces.
When $q=2$ the set of $n$--simplices is given by the set of pairwise commuting $n$--tuples of group elements.
As shown in \cite{O14,O15} for extraspecial $p$--groups $B(2,G)$ has non-trivial higher homotopy group.
In this paper we completely determine the homotopy type of the universal cover of $B(2,G)$ when $G$ is an extraspecial $p$--group. The space $B(2,G)$ is a classifying space for transitionally commutative principal $G$--bundles \cite[Theorem 2.2]{AG15}. Therefore our result implies the existence of non-trivial transitionally commutative principal $G$--bundles over the $r$--sphere $\sS^r$  where $2r$ is the rank of the Frattini quotient of $G$.  Another motivation is to understand the cohomology ring of $B(2,G)$ in order to gain a better understanding of the cohomology ring of the group $G$. 
The case of extraspecial $p$--groups is particularly interesting. Its cohomology ring in mod $p$ coefficients is completely determined by Quillen \cite{Q71a} when $p=2$. For odd $p$ only partial results are known, see \cite{BC92}. We believe that a description of $H^*(B(2,G),\FF_p)$ may help to gain more insight into the cohomology rings of extraspecial $p$--groups.   
 
Let $E(2,G)\rightarrow B(2,G)$ denote the pull-back of the universal principal $G$--bundle $EG\rightarrow BG$ along the inclusion $B(2,G)\subset BG$. It can be shown that $E(2,G)$ is homotopy equivalent to the nerve of the poset $\cC_G\aA(G)=\set{gA|\; g\in G \text{ and } A \subset G \text{ is an abelian subgroup}}$ ordered under the inclusion relation. In order to determine the homotopy type of the universal cover of $B(2,G)$ we study the coset poset  $\cC_G\aA(G)$.
The latter is weakly equivalent to a simpler coset poset $\cC_V\iI(V)$ where $V$ is a vector space over the finite field $\FF_p$ with a non-degenerate alternating bilinear form $\bi$, and $\iI(V)$ is the poset of isotopic subspaces of $V$.
Similar posets associated to a collection of subgroups   are studied extensively by various authors  \cite{S69,L74,Q78,V81,B00}. The techniques used in this paper span a variety of tools used in the study of topology of partially ordered sets  of subgroup collections. In particular there are two main ingredients, namely a combination of  the methods used by Vogtmann in \cite{V81} and by Quillen in \cite{Q78}.   Our main result is the following. 

\Thm{\label{X}Let $E$ be an extraspecial $p$--group. There is a fibration sequence
$$
\bigvee^{d(p,r)} \sS^r \rightarrow B(2,E) \rightarrow B\pi
$$
where $2r\geq 4$ is the rank of the Frattini quotient of $E$, 
$$
(-1)^rd(p,r)+1=(-1)^r p^{2r+1+r^2}+ \sum_{j=1}^r (-1)^{r-j} p^{2r+1-j+(r-j)^2} \left( \prod_{t=0}^{j-1} \frac{p^{2r-t} - p^{t}}{p^j-p^t} \right)
$$
and $\pi$ is the kernel of the multiplication map
$$
1\rightarrow \pi \rightarrow E\times E \rightarrow E/[E,E]\rightarrow 1.
$$
} 
\noindent There are two important consequences. The first one is the existence of non-trivial transitionally commutative $E$--bundles over $\sS^r$ as constructed in Theorem \ref{nontrivial_construction}. We also give an explicit construction of such a bundle.  
Another immediate application is a description of the cohomology of $B(2,E)$ in low degrees.

\Cor{\label{Y} Let $R$ denote a commutative ring, and $E$   an extraspecial $p$--group. Then there is an isomorphism
$$
H^i(B(2,E),R) \cong H^i( \pi,R)\;\; \text{ for } i<r
$$
where $2r\geq 4$ is the rank of the Frattini quotient of $E$ and $\pi$ denotes the fundamental group of $B(2,E)$.
}
\noindent We note that $\pi$ is completely determined in \cite{O14,O15}. It is a central extension of $E$ by a cyclic group of order $p$. For $p=2$ it splits as a direct product $E\times \ZZ/2$ hence its cohomology is known as a consequence of \cite{Q71a}.

The structure of the paper is as follows. 
We start with preliminary results on the homotopy theory of posets and introduce coset posets in \S \ref{preliminaries}.
In \S \ref{decomposingcoset} we describe the decomposition theorem (Theorem \ref{exact}) for coset posets. 
Next we specialize to the poset of abelian subgroups of extraspecial groups and prove some preliminary results in \S \ref{cosetabelian}.
The main result  of the paper is proved in \S \ref{cosetisotropic} in which we prove that the universal cover of $\cC_V\iI(V)$ is spherical when $G$ is an extraspecial $p$--group. Finally in \S \ref{secapp} we apply our results to $B(2,G)$.

\subsection*{Acknowledgement} The author would like to thank Alejandro Adem and Erg\" un Yal\c c\i n for helpful discussions. 

\section{Preliminaries} \label{preliminaries}
In this section we recall some results on the homotopy theory of posets. We are concerned with  special types of posets called coset posets. They arise naturally in the study of homotopy colimits of classifying spaces.

\subsection{Homotopy theory of posets} \label{homotopy_poset} Let $\xX$ denote a poset. 
  We will regard a poset  as a category with morphisms given by $x\rightarrow y$ whenever $x\leq y$. 
Notationally we usually do not distinguish the poset $\xX$ from its nerve $N(\xX)$ when talking about topological properties.  
  A poset map $f:\xX\rightarrow \yY$ can be regarded as a functor between the associated categories. 
A very useful notion associated to a map of posets is its fiber. 
The fiber of $f$ over an object $y\in \yY$  is defined by
$$
\underf{f}{y}=\set{x\in \xX|\; f(x)\leq y}, \text{ and dually }\; \overf{f}{y}=\set{x\in \xX|\; f(x)\geq y}.
$$
When interpreted as a category $\underf{f}{y}$ corresponds to the comma category $\comma{f}{y}$, and $\overf{f}{y}$ corresponds to $\comma{y}{f}$.
If $f$ is the identity map we simply write $\xX_{\geq x}$ for $\overf{f}{x}$. Other variations are $\xX_{>x}$, $\xX_{\leq x}$, and $\xX_{<x}$. Some constructions we will need are the following.
The join  $\xX\join\yY$ of the posets $\xX$ and $\yY$ is defined to be the poset whose underlying set is  the disjoint union $\xX\coprod \yY$. The ordering agrees with the given orderings on $\xX$ and $\yY$, and $x\leq y$ for every $x\in \xX$ and $y\in \yY$. 
The (unreduced) suspension $\Sigma \xX$ is defined to be the join $\set{0,1}\join \xX$ where $\set{0,1}$ has only the reflexive relation at each object.  

Next we recall some results on the homotopy theory of posets due to Quillen.
Let $f,g:\xX\rightarrow \yY$ be two  poset maps. If the maps satisfy the condition that $f(x)\leq g(x)$ (or dually $f(x)\geq g(x)$) for all $x\in \xX$ then they are homotopic, which will be  denoted by $f\simeq g$. This  property implies that a poset with an initial (or terminal) object is contractible since the identity map will be homotopic to the constant map.  Another  fact  we will use is a version of Quillen's Theorem A for posets.

\Thm{\cite[Proposition 1.6]{Q78}\label{QuillenThA} Let $f:\xX\rightarrow \yY$ be a map of posets. If  the fiber $\overf{f}{y}$ (or $\underf{f}{y}$)   is contractible for all objects $y$ in $\yY$ then $f$ is a homotopy equivalence. }

The dimension $\dim (\xX)$ of a poset $\xX$ is defined to be the supremum of the integers $k$ such that there is a chain $x_0< x_1 < \cdots < x_k$. 
A poset of dimension $n$ is said to be   $n$--spherical if it (or rather its nerve) is $(n-1)$--connected. 
In the study of spherical posets the following theorem of Quillen is very useful. This version is the dual of the one stated in \cite[Theorem 9.1]{Q78}.

\Thm{\cite[Theorem 9.1]{Q78} \label{Quillen} Let $f:\xX\rightarrow \yY$ be a map of posets.
If $\yY$ is $(\dim \yY)$--spherical, $\overf{f}{y}$ is $\dim (\yY_{\geq y})$--spherical and  $\yY_{<y}$ is $\dim (\yY_{<y})$--spherical for all $y$ then $\xX$ is $(\dim\yY)$--spherical.
}

\subsection{Coset posets}\label{coset_poset}  Let $G$ be a finite group and $\fF$   a collection of subgroups  of $G$. We define the coset poset associated to the collection $\fF$ as a set
$$
\cC_G\fF=\set{gA|\;A\in \fF \text{ and } g\in G}
$$
and regard it as a poset ordered under inclusion.
 There is a morphism $gA\rightarrow g'A'$ whenever $gA$ is contained in $g'A'$ as a set.
 Note that $\cC_G\fF$ is a poset with $G$--action given by left translation and its nerve is a $G$--space i.e.  a simplicial set with $G$--action. 
We also define a relative version. Let $H\subset G$ be a subgroup. 
Given a (left) coset $gH\subset G$ we define  
$$
\cC_{gH} \fF=\set{ghA|\; A\in \fF\text{ and }gh\in gH}.
$$ 
The poset structure is similarly induced by inclusion as sets, and  $gHg^{-1}$ acts by left translation. 
An element $x\in G$ induces a map $x:\cC_{gH}\fF\rightarrow \cC_{xgH} \fF$ defined by $ghA\mapsto xghA$.

We describe a basic result which usually allows us to work with smaller collections of subgroups.   
Let $\fF_{\cap H}$ denote the collection $\set{A\cap H| A\in \fF}$.  
There is  an $H$--equivariant map
$$
i_H:\cC_H \fF \rightarrow \cC_{H} \fF_{\cap H}
$$
defined by $hA\mapsto h(A\cap H)$. Under suitable conditions this map turns out to be a weak equivalence. In fact, it becomes a weak $H$-equivalence. This means that the induced map on fixed points is a weak equivalence for all subgroups.  
\Def{\rm{We say $\fF$ is $H$--\textit{stable} if for all $A\in \fF$ the intersection $A\cap H $ is also in $\fF$.  }}
\noindent For example the collection of abelian subgroups satisfies this property with respect to any subgroup. 
\Pro{\label{H-eq} If $\fF$ is $H$--stable then  $i_{H}:\cC_H\fF\rightarrow \cC_H \fF_{\cap H}$ is a weak $H$--equivalence.}
\Proof{Note that by assumption we have $\fF_{\cap H}\subset \fF$.
Let $K\subset H$ and consider the restriction of $i_H$ to the fixed points $i_H:(\cC_H\fF)^K\rightarrow (\cC_H \fF_{\cap H})^K$. 
If the latter fixed point set   is not empty then  the fiber over a coset $hB\in \cC_H \fF_{\cap H}$ fixed under the action of $K$   is given by
\begin{align*}
\overf{(i_H)}{hB} &=\set{h'A\in (\cC_H\fF)^K|\; h'(A\cap H)  \supset hB} \\
&= \set{hA\in (\cC_H\fF)^K|\; A\cap H  \supset B}\\
&= \set{hA\in \cC_H\fF|\; A\cap H  \supset B \text{ and }  K\subset hAh^{-1}}.
\end{align*}
Note that if $h'(A\cap H) \supset hB$ then $h^{-1}h' (A\cap H)\supset B$. In particular, we have $1\in h^{-1}h' (A\cap H)$, that is, $h'=hx$ for some $x\in A\cap H$. Therefore we have $h'A =hA$, and the first equality above follows.
The fiber is contractible since $hB$ is initial  and the result follows from Theorem \ref{QuillenThA}.  
}

\subsection{Homotopy colimit of classifying spaces} \label{hocolim_classifying} 
Our reference for homotopy colimits of classifying spaces is \cite[\S 3]{O14}.
Let $\catS$ denote the category of simplicial sets.
The classifying space of a group  can be seen as a functor $B:\catGrp \rightarrow \catS$ from the category of groups to the category of simplicial sets. Let $\fF$ denote a collection of subgroups of $G$. We regard $\fF$ as a partially ordered set. Restricting the classifying space functor gives $B:\fF \rightarrow \catS$. There is a fibration sequence
\begin{equation}\label{fibseq_G}
\hocolim{\fF} (G/-) \rightarrow \hocolim{\fF} B \rightarrow BG 
\end{equation}
induced by the inclusions $BA\subset BG$ for each $A\in \fF$, where $(G/-):\fF\rightarrow \catS$ is the functor which sends $A$ to the coset $G/A$ regarded as a discrete simplicial set. The fiber can be described equivalently as follows. The homotopy colimit of the cosets $G/A$ as $A$ runs over the poset $\fF$ is isomorphic to the nerve of the coset poset
$$
\cC_G \fF =\set{\,gA|\; g\in G \text{ and } A\in \fF\,}
$$
(see \cite[Proposition 5.12]{DH01}).  
The homotopy long exact sequence of the fiber sequence \ref{fibseq_G} gives an exact sequence
$$
1\rightarrow \pi_1\, \cC_G \fF  \rightarrow \pi_1\, \hocolim{\fF} B \rightarrow G \rightarrow 1
$$
of fundamental groups.

\Pro{\label{far}
Assume that the collection $\fF$ has an initial object. Then
$$
\pi_1\, \hocolim{\fF} B \cong \colim{A\in\fF} A
$$
where the colimit is taken in the category of groups.
}
\Proof{This is a consequence of \cite[Corollary 5.1]{Far04}.}

Let $\pi$ denote the colimit of the groups $A$ in $\fF$. 
There is a commutative diagram of groups
$$
\begin{tikzcd}
\pi  \arrow{r}{\pi} & G \\
A \arrow{u}{i_A} \arrow{ur}  &
\end{tikzcd}
$$
which implies that the natural map $i_A$ is a monomorphism. Therefore we can regard subgroups in $\fF$ as subgroups of $\pi$, and talk about the cosets $\pi/A$. Consider the fibration sequence
$$
\hocolim{\fF} (\pi/-) \rightarrow \hocolim{\fF} B \rightarrow B\pi
$$
induced by the inclusions $BA\subset B\pi$. Note that the fiber is a simply connected space, and it can be identified with the nerve of the coset poset
$$
\cC_\pi \fF=\set{\,gA|\; g\in \pi \text{ and } A\in \fF\,}.
$$
We will study this object using homotopy theoretic methods for posets.

\section{Decomposing coset posets} \label{decomposingcoset}

We will describe a decomposition of $\cC_G\fF$ as a homotopy colimit. For basic properties of homotopy colimits we refer to \cite{Dwy97, BK72,GJ99}. 

Let $f:\xX\rightarrow \yY$ be a map of posets. Define a functor $\overf{f}{-}:\yY^\op\rightarrow \catS$ by sending an object $y$ to the nerve of the fiber $\overf{f}{y}$. Here $(-)^\op$ denotes the opposite category.   There is a natural map
\begin{equation}\label{nat}
\wtilde f :\hocolim{ } \overf{f}{-} \rightarrow \xX
\end{equation}
induced by the inclusions $\overf{f}{y}\rightarrow\xX$. This map is a weak equivalence  \cite[\rom{4} \S 5.1]{GJ99}. 
We need some notation for the next result. Let $\lL(G)$ denote the collection of all subgroups including the trivial subgroup and the group $G$. For a normal subgroup $H$ we define $\fF^{\vee H}$ to be the collection $\set{A\in \fF| AH=G}$ where $AH\subset G$ denotes the subgroup generated by $A$ and $H$. 
Let $xH$ denote a coset in $G$. There is a map of posets
\begin{equation}\label{thetaH}
\theta_{xH}: \cC_{G} \fF^{\vee H} \rightarrow \cC_{xH} \fF
\end{equation}
defined by sending a coset $yA$  to the coset $xhA$
where  $y=xha$ for some $h\in H$ and $a\in A$. 

\begin{pro}\label{cyclic} Let $1\rightarrow H\rightarrow G \stackrel{\pi}{\rightarrow} \bar G \rightarrow 1$ be an exact sequence of groups where $\bar G$  is isomorphic to a cyclic group $C_p$ of prime order $p$. Choose an element $g\in G$ such that $\bar{g}=\pi(g)$ generates $\bar{G}$. Assume that $\fF$ contains the trivial subgroup and a subgroup not contained in $H$.
Then there is a weak equivalence
$$
\wtilde\pi:\hocolim{\cC_{C_p}\lL(C_p)^\op}{\overf{\pi}{-}} \rightarrow \cC_G\fF 
$$ 
and  the functor $\overf{\pi}{-}$ can be identified as
\begin{equation}\label{fibers}
\overf{\pi}{\bar G}=\cC_G\fF^{\vee H}\;\;\text{ and } \;\overf{\pi}{\bar g^t}=\cC_{g^tH}\fF.
\end{equation}
\end{pro}
\begin{proof}
The map $\wtilde \pi$ is obtained as follows.
Let $\pi\fF$ denote the poset $\set{\pi(A)|\; A\in \fF}$ of subgroups of $\bar{G}$. There is an induced map of  posets $\pi:\cC_G\fF\rightarrow \cC_{\bar G}(\pi\fF)$  defined by $gA\mapsto \pi(gA)$. The decomposition \ref{nat} applied to $\pi$ gives a weak equivalence
\begin{equation}\label{decomposition}
\wtilde \pi: \hocolim{} \overf{\pi}{-} \rightarrow \cC_G\fF
\end{equation}
where the colimit is over the poset $\cC_{C_p} \lL(C_p)^{\op}$. By the assumptions on $\fF$ this poset can be identified with  
$$\cC_{\bar G}(\pi\fF)^\op=\set{\bar g, \bar g^2,\cdots, \bar g^p ,\bar G}$$
 partially ordered under reverse inclusions. There are two types of fibers $\overf{\pi}{\bar G}$ and $\overf{\pi}{\bar g^t}$ where $1\leq t\leq p$.   Then $\cC_G\fF$ is the homotopy colimit of the diagram
\begin{equation}\label{picture}
\begin{tikzcd}
\cC_{gH}\fF & \cdots &  \cC_{g^pH}\fF \\
 &\arrow{ul}{\theta_{gH}} \cC_G\fF^{\vee H}\arrow[ur,"\theta_{g^pH}"']&
\end{tikzcd}
\end{equation}
where $\theta_{g^tH}$ for $t=1,2,\cdots,p$ denote the natural inclusions defined above. 
\end{proof} 
 
\Ex{\rm{Let $V$ be a vector space over $\FF_p$, and let $\tT(V)$ denote the poset of proper subspaces. Choose a subspace $W$ of codimension one. Since $\tT(V)$ is $W$-stable 
$$
i_W: \cC_W\tT(V) \rightarrow \cC_W\tT(V)_{\cap W}. 
$$
is a weak equivalence by Proposition \ref{H-eq}.
Moreover $\tT(V)_{\cap W}$ is contractible since  $W$ is  terminal. Therefore the  decomposition in Proposition \ref{cyclic}  for $W = H$ gives
$$
\cC_V \tT(V) \simeq \bigvee^{p-1}\Sigma (\cC_V\tT(V)^{\vee W}).
$$ 
}
}

\subsection{Homotopy sections} 
We make an   assumption in addition to the set-up in Proposition \ref{cyclic}. Let $g^{t} \fF^{\vee H}$ denote the poset of left translations $\set{g^tA|\;A\in \fF^{\vee H}}$ regarded as a subposet of $\cC_G\fF^{\vee H}$.
Assume that there exists a map of posets
\begin{equation}\label{s}
s:\cC_H\fF \rightarrow \fF^{\vee H}
\end{equation}
such that $\theta_H s \simeq \ast$  and the composition 
$$
\cC_H \fF \stackrel{s}{\longrightarrow} \fF^{\vee H} \stackrel{g^{-1}}{\longrightarrow} g^{-1} \fF^{\vee H} \stackrel{\theta_H}{\longrightarrow} \cC_H \fF
$$
is homotopic to the identity map. Here $g^{-1}:\fF^{\vee H} \rightarrow g^{-1}\fF^{\vee H}$ is defined by $A\mapsto g^{-1}A$, and $\theta_H$ is restricted to $g^{-1} \fF^{\vee H}$. The restricted map sends $g^{-1}A\in g^{-1} \fF^{\vee H}$ to $hA$ where $g^{-1}=ha$ for some $h\in H$.
 
\Lem{\label{hsec_exists}Given the map $s$ in \ref{s} there exists a homotopy section
$$
s_{g^kH}  :\cC_{g^kH}\fF \rightarrow   \cC_G\fF^{\vee H}
$$  
of $\theta_{g^kH}$ such that $\theta_{g^{k-1}H}s_{g^kH}\simeq \ast$  for all $1\leq k\leq p$ . 
}  
\begin{proof}
We define $s_{g^kH} $ to be the composite
$$
s_{g^kH}  :\cC_{g^kH}\fF \stackrel{g^{-k}}{\longrightarrow} \cC_{H}\fF 
\stackrel{s}{\longrightarrow}  \fF^{\vee H} \stackrel{g^{k-1}}{\longrightarrow}  g^{k-1}\fF^{\vee H} \subset \cC_G\fF^{\vee H}
$$  
for $1\leq k\leq p$. Then we have
\begin{equation}\label{section_property}
\theta_{g^lH}s_{g^kH}= g^l \theta_H g^{k-1-l}sg^{-k} \simeq \left\lbrace\begin{array}{lll}
\idy & \text{ if } l=k\\
\ast & \text{ if } l=k-1.
\end{array} \right.
\end{equation}
where we have used $\theta_{g^lH}=g^l \theta_H g^{-l}$.
This gives a sequence of sections
$$
\begin{tikzcd}
     &     &\arrow{dl}{\theta_{g^{-1}H}}  \cC_{G} \fF^{\vee H} \arrow[dr, "\theta_{H}"'] &    & \arrow{dl}{\theta_{H}} \cC_{G} \fF^{\vee H} \arrow[rd, "\theta_{gH}"']& &    \\
      & \cdots \cC_{g^{-1}H}\fF   &                 & \cC_{H}\fF   \arrow[lu, dotted, bend right,"s_{H}"'] &                        &  \cC_{gH}\fF  \arrow[lu, dotted, bend right, "s_{gH}"']  \cdots &
\end{tikzcd}
$$
such that $\theta_{g^{k-1}H} s_{g^kH} \simeq 0$.
\end{proof}

For a poset $\xX$ we write $H_*(\xX)$ (and $\tilde H_*(\xX)$) for the (reduced) integral homology groups of the nerve of $\xX$. Given a map of posets $f:\xX\rightarrow \yY$ the induced map in homology is denoted by $f_*: H_*(\xX)\rightarrow H_*(\yY)$ (similarly for   reduced homology groups).

\Thm{\label{exact}
Let $1\rightarrow H\rightarrow G \stackrel{\pi}{\rightarrow} \bar G \rightarrow 1$ be an exact sequence of groups where $\bar G$  is isomorphic to a cyclic group $C_p$ of prime order $p$ generated by $\bar g$. Let  $\fF$ be a collection of subgroups containing the trivial subgroup and a subgroup not contained in $H$. Assume that there is a map $s:\cC_H\fF\rightarrow \fF^{\vee H}$  such that $\theta_Hg^{-1}s\simeq \idy$ and $\theta_H s\simeq \ast $. Then there is a split exact sequence
$$
0\rightarrow \tilde H_i(\cC_G\fF) \rightarrow \tilde H_{i-1}(\cC_G\fF^{\vee H})^{\oplus(p-1)} \stackrel{\theta}{\rightarrow} \bigoplus_{k=1}^p \tilde H_{i-1}(\cC_{g^kH} \fF) \rightarrow 0,\;\;i\geq 1 
$$
where $\theta(x_t)=(\theta_{g^{t-1}H})_*(x_t)-(\theta_{g^tH})_*(x_t)$ for $1\leq t< p$.}
 
\noindent This theorem is an immediate consequence of Lemma \ref{hsec_exists}, the homotopy colimit decomposition \ref{picture}, and the following result based on the Mayer-Vietoris sequence. 
 
\Pro{\label{exact_space}
Let $X_0,X_1,\cdots,X_{m-1}$ denote a sequence of subspaces of a space $X$ such that their union gives $X$ and any pairwise intersection $X_k\cap X_l$ with  $k\not=l$ gives the same subspace denoted by $X_\text{int}$.  Assume that each inclusion $\iota_k: X_\text{int}\rightarrow X_k$ has a homotopy section $s_k$ such that  $\iota_{k-1}s_k\simeq \ast$  where $k$ is written mod $m$. Then the associated Mayer-Vietoris sequence splits to yield short exact sequences
$$
0\rightarrow \tilde H_{i}(X) \rightarrow \tilde H_{i-1}(X_\text{int})^{\oplus(m-1)} \stackrel{\theta}{\rightarrow} \bigoplus_{k=0}^{m-1} \tilde H_{i-1}(X_k) \rightarrow 0
$$
} 
\begin{proof}
There is  a Mayer--Vietoris sequence
$$
\cdots \rightarrow \tilde H_{i-1}(X_\text{int})^{\oplus m-1}   \stackrel{\theta}{\rightarrow} \bigoplus_{k} \tilde H_{i-1}(X_k)\rightarrow \tilde  H_{i-1}(X)\rightarrow \tilde  H_{i-2}(X_\text{int})^{\oplus m-1}  \rightarrow \cdots
$$
of reduced homology groups associated to the diagram of inclusions
$$
\begin{tikzcd}
X_0 & \cdots &  X_{m-1} \\
 &\arrow{ul} X_\text{int} \arrow{ur}&
\end{tikzcd}
$$
One can obtain this sequence as a special case of the Bousfield-Kan spectral sequence \cite[\S XII.5.7]{BK72}.
The map $\theta$ sends an element $x$ in the $k$-th summand to $(\iota_{k-1})_*(x)-(\iota_{k})_*(x)$.
We will show that the Mayer-Vietoris long exact sequence splits by constructing a section of $\theta$.
The sections $s_k$ can be put together to define a section of $\theta$. Let $\Delta:H_{i-1}(X_\text{int})\rightarrow H_{i-1}(X_\text{int})^{\oplus m-1}$ denote the diagonal map. Then  define
$$
s_\theta:  \bigoplus_{k=0}^{m-1} \tilde H_{i-1}(X_k) \rightarrow \tilde H_{i-1}(X_\text{int})^{\oplus m-1} 
$$ 
by $s_\theta(\alpha_k)=(s_{k})_*(\alpha_k)$ if $k>0$ and $s_\theta(\alpha_0)=\Delta (s_{0})_*(\alpha_0)$. Then $s_\theta$ is a section of $\theta$. 
\end{proof}

\Cor{\label{exact-n} Assume that the map $s$ exists as defined in Theorem \ref{exact}. If $\cC_G\fF^{\vee H}$ and $\cC_H \fF$ are $(n-1)$--spherical then $\cC_G\fF$ is $n$--spherical.}
\Proof{We set $X_k=\cC_{g^kH} \fF$ and use notation from Proposition \ref{exact_space}.
Van Kampen's theorem gives an isomorphism
\begin{equation}\label{vankampen}
\pi_1(X) \cong (\ast_{k} \pi_1(X_k))/\sim
\end{equation} 
where the quotient relation is generated by $(\iota_k)_*(x)\sim(\iota_{k-1})_*(x)$ for $x\in \pi_1(X_\text{int})$ and $0\leq k\leq m-1$. We prove now that $X$ is simply connected by a similar argument to that used in the proof of Proposition \ref{exact_space}. Fix $k$ and let $x_k$ be an element of $\pi_1(X_k)$. Taking $x=(s_{k})_*(x_k)$ we see that $x_k=(\iota_{k})_*(x)\sim (\iota_{k-1})_*(x)=1$.   Therefore the equivalence relation in \ref{vankampen} identifies every generator with the identity element, and as a result we have $\pi_1(X)=1$. Combining this observation with Theorem \ref{exact} gives the desired result.
}

\section{Coset poset of abelian subgroups of extraspecial groups} \label{cosetabelian}

Let $\aA(G)$ denote the poset of abelian subgroups of $G$. We will study $\cC_G\aA(G)$ when $G$ is an extraspecial $p$--group. The universal cover in this case turns out to be spherical. 

\subsection{Extraspecial groups} An extraspecial $p$--group is a central extension of the form
$$
0\rightarrow \ZZ/p \rightarrow E \stackrel{\nu}{\rightarrow} V \rightarrow 0
$$
where $V$ is an elementary abelian $p$--group of rank $2r$. The kernel of the extension is both the center  and the commutator  of the group \cite[Chapter 8]{As86}. We will denote the center by $Z(E)$. The quotient $V$ is also the Frattini quotient. The commutator induces a bilinear form on $V$: $\bi(v_1,v_2)=[\nu^{-1}(v_1),\nu^{-1} (v_2)]$.  
A subspace $I\subset V$ is called isotropic if the restriction of $\bi$ on $I$ is zero. Let $\iI(V)$ denote the collection of isotropic subspaces of $V$. 
There is a map of posets $\aA(E)\rightarrow \iI(V)$ defined by $A\mapsto \nu(A)$ which induces a map of posets  between the coset posets $\hat\nu:\cC_E\aA(E)\rightarrow \cC_V\iI(V)$ defined by $xA\mapsto \nu(xA)$. 

\Pro{\label{reduction}The induced map $\hat\nu:\cC_E\aA(E)\rightarrow \cC_V\iI(V)$ is a weak equivalence.}
\Proof{The fiber $\underf{\hat \nu}{v+I}$  contains the coset $\nu^{-1}(v+I)$ as a terminal object. Hence the given map is a weak equivalence by Theorem \ref{QuillenThA}.}

\subsection{$r=1$ case}\label{r1} Consider the coset poset $\cC_V\iI(V)$ when $r=1$. Then $\iI(V)$ consists of all the one-dimensional subspaces and the zero subspace. The corresponding space is a one-dimensional connected space hence has the homotopy type of a wedge of circles. For example, the set of abelian subgroups for the quaternion group $Q_8$ is given by $\set{\Span{i},\Span{j},\Span{k},\Span{-1},1}$. Under the quotient map $\nu:E\rightarrow V$ the maximal abelian subgroups map to the one-dimensional subspaces, the center and the identity subgroup both map to the zero subspace. In effect the weak equivalence in Proposition \ref{reduction} implies that to determine the homotopy type of the coset poset we can restrict to maximal abelians and their intersections.

\subsection{Fundamental group} In \cite[\S 4.2]{O15} the colimit of abelian subgroups of an extraspecial $p$--group is computed. Recall from section \ref{hocolim_classifying} that this colimit is isomorphic to the fundamental group of the homotopy colimit of the classifying spaces $BA$ of abelian subgroups $A\subset E$. 

\Pro{\label{colim_ext} \cite{O15} Assume $r\geq 2$. There is an isomorphism of groups
$$
\phi:\colim{A\in\aA(E)} A \rightarrow  \pi 
$$
where $\pi$ is the kernel of the multiplication map $m(e,e')=ee'[E,E]$: 
$$
1\rightarrow \pi \rightarrow E\times E \stackrel{m}{\longrightarrow} E/[E,E] \rightarrow 1
$$
and $\phi(a)=(a,a^{-1})$ for $a\in A$.
} 

Note that $\pi$ is an extension of $E$ by a cyclic group of order $p$. Using the isomorphism $\phi$, we see that the homotopy long exact sequence of \ref{fibseq_G} gives
$$
1\rightarrow \pi_1\, \cC_E\aA(E) \rightarrow \pi \stackrel{\pi_1}{\longrightarrow} E \rightarrow 1
$$
where $\pi_1$ is the projection onto the first factor.
Therefore we obtain the following.

\Cor{\label{pi1} Assume $r\geq 2$.
There is an isomorphism
$$
\pi_1\, \cC_E \aA(E) \cong \ZZ/p.
$$
}

\subsection{Heisenberg group}\label{Heisenberg} Let $H(V)$ denote the Heisenberg group associated to the vector space $V$ with the bilinear form $\bi$. This group is defined to be  the set $ V\times \ZZ/p$ with the multiplication rule given by
$$
(v_1,t_1)(v_2,t_2)=(v_1+v_2,\bi(v_1,v_2)+t_1+t_2).
$$
Let us identify $Z(E)$ with $\ZZ/p$ via an isomorphism. 

\Lem{\label{pi_HV} 
There is a surjective homomorphism of groups defined by 
$$
\varphi:\pi \rightarrow H(V), \;\;\;\; (e,e')\mapsto (\nu(e),ee')
$$
whose kernel consists of pairs $(a,a^{-1})$ where $a\in Z(E)$.
}
\Proof{
We check that $\varphi$ is a group homomorphism. On one hand we have
$$
\varphi((e_1,e_1')(e_2,e_2'))=\varphi(e_1e_2,e_1'e_2')= (\nu(e_1e_2),e_1e_2e_1'e_2')
$$
and on the other hand
$$
\varphi((e_1,e_1'))\varphi((e_2,e_2'))=(\nu(e_1),e_1e_1')(\nu(e_2),e_2e_2')=(\nu(e_1)\nu(e_2),[e_1,e_2]e_1e_1'e_2e_2').
$$
Note that both give the same result since 
$$[e_1,e_2]e_1e_1'e_2e_2' = e_1e_2e_1^{-1}e_2^{-1}e_1e_1'e_2e_2'= e_1e_2 e_1'e_2^{-1}e_2e_2'=e_1e_2 e_1'e_2'$$
where we have used the fact that $e_1e_1'$   belongs to the center $Z(E)$. To see that it is surjective observe that the images of the elements of the form $(e,e^{-1})$ generate $H(V)$. The kernel is as described by definition of the map.
}

Consider the projection $\pr_1:H(V)\rightarrow V$ onto the first factor. An isotropic subspace $I$ can be identified with an abelian subgroup of $H(V)$ via  the map $a\mapsto (a,0)$. We can regard $\iI(V)$ as a collection of subgroups of $H(V)$. First we need a preliminary result. Regard the classifying space as a functor $B$ restricted on the posets $\aA(E)$ and $\iI(V)$. We will compute the fundamental group of the homotopy colimit of  $B:\iI(V)\rightarrow \catS$.

\Pro{\label{fund_hoco_iso}
There is a natural isomorphism
$$
\pi_1 \hocolim{\iI(V)} B \cong H(V).
$$
}
\begin{proof}
 The natural map of posets $\aA(E)\rightarrow \iI(V)$ given by $A\mapsto \nu(A)$ induces a map
$$
\hocolim{\aA(E)} B \rightarrow \hocolim{\iI(V)} B
$$
between the homotopy colimits. The homotopy fiber of this map is $B\ZZ/p$. This can be seen from the diagram
$$
\begin{tikzcd}\label{diag_fib}
                           & B\ZZ/p \arrow{d} \arrow[r,equal] & B\ZZ/p \arrow{d} \\   
\cC_E\aA(E)\arrow{r}\arrow{d}{\simeq} & \hocolim{\aA(E)}B \arrow{r} \arrow{d} & BE \arrow{d}\\
\cC_V\iI(V)\arrow{r}  & \hocolim{\iI(V)}B \arrow{r}  & BV  
\end{tikzcd}
$$
where we have used Proposition \ref{reduction} to identify the horizontal fibers found in \ref{fibseq_G}.
Both homotopy colimits have the same universal cover and they differ by $B\ZZ/p$. Therefore all the sequences in the diagram are fibration sequences.
Under the identifications of Proposition \ref{colim_ext} and Corollary \ref{pi1} consider the diagram of fundamental groups
$$
\begin{tikzcd}
& & 0\arrow{d} & 0\arrow{d} & \\
& & \ZZ/p \arrow[r,equal] \arrow{d}& \ZZ/p \arrow{d} & \\
0\arrow{r} &\ZZ/p\arrow[r]\arrow[d,equal] & \pi \arrow{r}{\pi_1} \arrow{d}& E \arrow{d}\arrow{r} &0\\
0\arrow{r} & \ZZ/p\arrow[r]  & \bar \pi \arrow{r} \arrow{d}   & V  \arrow{r}\arrow{d}&0\\
&& 0 & 0 
\end{tikzcd}
$$
where $\bar \pi$ is the fundamental group of the homotopy colimit of $B:\iI(V)\rightarrow \catS$.
The middle column is the extension corresponding to the homomorphism in Lemma \ref{pi_HV} since the cyclic subgroup $\ZZ/p$ in $\pi$ is
precisely $\ker\varphi$. Therefore $\bar \pi$ is isomorphic to $H(V)$ as desired.

\end{proof}

As an immediate consequence of this computation we have the following.

\Cor{\label{covering}  Assume that $r\geq 2$ where $2r$ is the dimension of $V$. Then the natural map 
$$\hat\pr_1:\cC_{H(V)}\iI(V)\rightarrow \cC_V\iI(V)$$
 induced by the projection $\pr_1:H(V)\rightarrow V$ is the universal covering map.
} 
\begin{proof}

The projection map $\pr_1:H(V)\rightarrow V$ induces a map of fibrations
$$
\begin{tikzcd}
\cC_{H(V)} \iI(V) \arrow{r}\arrow{d} & \hocolim{\iI(V)}B \arrow{r} \arrow[d,equal] & B (H(V)) \arrow{d}\\
\cC_{V} \iI(V) \arrow{r} & \hocolim{\iI(V)}B \arrow{r}  & B V 
\end{tikzcd}
$$
from which we can conclude that the map between the fibers is the universal covering map.
\end{proof}

Next we summarize the relationship between the coset posets associated to $\aA(E)$ and $\iI(V)$.

\begin{pro}\label{diag_covers}
The homomorphism $\varphi: \pi \rightarrow H(V)$ induces a diagram 
$$
\begin{tikzcd}
\cC_\pi \aA(E) \arrow{r} \arrow{d}{\simeq} & \cC_E \aA(E) \arrow{d}{\simeq} \\
\cC_{H(V)} \iI(V) \arrow{r} & \cC_V \iI(V)
\end{tikzcd}
$$
between the universal covering maps whose fibers are given by a copy of $\ZZ/p$.
\end{pro}

Our aim is to prove the following result which will follow from the corresponding result (Theorem \ref{main1}) for the coset poset $\cC_{H(V)}\iI(V)$ in view of Proposition \ref{diag_covers}. 
 
\Thm{\label{main0} Assume $r\geq 2$. The coset poset $\cC_{\pi} \aA(E)$ is $r$--spherical.
} 

We leave the proof of this theorem to Section \ref{cosetisotropic} and conclude this section with a formula of the Euler characteristic involving dimensions of Steinberg representations of symplectic groups.

\subsection{Euler characteristic} 
The geometric realization of the poset  $\iI^\circ(V)=\iI(V)-\set{0}$ can be identified with the Tits building of the symplectic group $\Sp(V)$. It is $(r-1)$--spherical by the Solomon-Tits theorem \cite{S69}  where $2r=\dim(V)$.
The top dimensional homology affords the Steinberg representation of $\Sp(V)$. We will describe the top dimensional homology of $\cC_{H(V)} \iI(V)$ as the kernel of a long exact sequence   which consists of Steinberg representations for symplectic groups.

As in Section \ref{hocolim_classifying} the coset space $\cC_{H(V)} \iI(V)$ can be identified with the  homotopy colimit of the functor $HV/-:\iI(V)\rightarrow \catS$ which sends an isotropic subspace $I$ to the coset $H(V)/I$.
There is a Bousfield-Kan spectral sequence \cite[\S XII.5.7]{BK72} computing the integral homology of this homotopy colimit whose $E_2$--page is given by
 $$
\text{colim}_i H_j(HV/I) \Rightarrow H_{i+j}(\cC_{H(V)} \iI(V))
$$
which collapses onto the $j=0$ axis since coset spaces are discrete. Therefore we have
$$
 H_{k}(\cC_{H(V)} \iI(V)) \cong \text{colim}_k\, \ZZ[HV/-]
$$ 
where $\ZZ[HV/-]$ is the functor $\iI(V)\rightarrow \catAb$ which sends $I$ to the free abelian group generated by the coset $HV/I$. We will describe a spectral sequence to compute such derived colimits.  For details of (the dual version of) this construction see \cite[\S 4]{O14}

The derived colimits of a functor $F:\pP \rightarrow \catAb$ can be described as the homology $H_*(\pP,F)$ of a chain complex  
$$C_k(\pP;F)= \bigoplus_{p_0<p_1<\cdots <p_k} F(p_0) $$
whose differential is induced by removing an object $p_0<p_1<\cdots < \hat p_i<\cdots <p_k$.
We  describe a nice filtration on the poset.   For an object $p\in \pP$ we define   $\height(p)=-\dim(\pP_{\geq p})$. This yields a filtration on the chains 
$$
C_*(\pP;F_r) \subset \cdots \subset C_*(\pP;F_1)\subset C_*(\pP;F_0)
$$
where the functors $F_i$ are  defined as follows: $F_i(p)=0$ if $\height(p)>-i$, and $F_i(p)=F(p)$ otherwise. 
Associated to the filtration there is a spectral sequence in the eighth octant:
$$
E_{i,j}^1=\bigoplus_{p\in \pP|\height(p)=-i}H_{i+j}((\pP_{\geq p},\pP_{>p}),F(p)) \Rightarrow H_{i+j}(\pP,F).
$$ 
Now we apply this spectral sequence to the poset $\pP=\iI(V)$ and the functor $F=\ZZ[HV/-]$.  
The first page consists of a direct sum of homology groups
\begin{equation}\label{E1}
H_{i+j}((\iI(V)_{\geq I},\iI(V)_{>I}),F(I))\cong  \ZZ [HV/I]\otimes H_{i+j}(\Sigma\,\iI^\circ(I^\perp/I))
\end{equation}
over isotropic subspaces $I$ of height equal to $-i$. Note that $\height(I)=\dim(I)-r$ hence in fact the sum is over $I$ of dimension $r-i$. In this case the quotient $I^\perp/I$ is of dimension $2i$ and  $\iI^\circ(I^\perp/I)$ is $(i-1)$--spherical. Therefore the term in \ref{E1} is only non-zero if $j=0$ and the spectral sequence collapses to give a long exact sequence
\begin{eqnarray*}\label{chains}
0\rightarrow   \ZZ[HV]\otimes H_{r-1}(\iI^\circ(V)) \rightarrow \bigoplus_{I_1} \ZZ[V/I_1]\otimes H_{r-2}(\iI^\circ(I_1^\perp/I_1))\rightarrow \cdots \\
 \rightarrow \bigoplus_{I_{r-1}} \ZZ[V/I_{r-1}]\otimes H_0(\iI^\circ(I_{r-1}^\perp/I_{r-1}))\rightarrow \bigoplus_{I_r} \ZZ[V/I_r] \rightarrow 0
\end{eqnarray*}
whose homology computes the homology of $\cC_{H(V)} \iI(V)$. Here the direct sum runs over isotropic subspaces $I_k$ of dimension $k$.
 
\Cor{\label{euler}The Euler characteristic of $\cC_{H(V)} \iI(V)$ is given by 
$$(-1)^rd(p,r)+1=(-1)^rp^{2r+1+r^2}+ \sum_{j=1}^r (-1)^{r-j} p^{2r+1-j+(r-j)^2} \left( \prod_{t=0}^{j-1} \frac{p^{2r-t} - p^{t}}{p^j-p^t} \right) .$$
}
\Proof{
Let $N_{j}$ denote the number of $j$-dimensional isotropic subspaces $I_j$ in $V$, and $D_j$   the dimension of the Tits building associated to the quotient $I_j^\perp/I_j$. When $j=r$ the quotient is the zero space and we set $D_r=1$. Then the Euler characteristic is given by the alternating sum
$$
\sum_{j=0}^r (-1)^{r-j} |HV/I_j| N_{j} D_{j}
$$
where the number of isotropic subspaces is given by
$$
N_{j}=\prod_{t=0}^{j-1} \frac{p^{2r-t} - p^{t}}{p^j-p^t}  
$$
and we take $N_j=1$ if $j=0$.   Note that the number of $j$ linearly independent isotropic vectors is given by $(p^{2r}-1)(p^{2r-1}-p)\cdots (p^{2r-(j-1)}-p^{j-1})$ since the dimension of the orthogonal complement of an isotropic subspace $I$ is equal to $2r-\dim(I)$. Dividing this product by the number of $j$ linearly independent vectors in a $j$--dimensional vector space gives the formula for $N_j$.
The dimension $D_{j}$ of the Tits building of $I_j^\perp/I_j$ equals the order of the unipotent radical  of the Borel subgroup of $\Sp(I_j^\perp/I_j)$ which is $p^{(r-j)^2}$ (see \cite[\S 3.5.4.]{Wil09}).  
}

\section{Coset poset of isotropic subspaces} \label{cosetisotropic}

In this section we focus on the coset poset $\cC_{H(V)} \iI(V)$ as introduced in \S \ref{cosetabelian}. In Theorem \ref{main1} we prove that this poset is $r$--spherical where $2r=\dim(V)$. We start with some results about subspace collections of vector spaces of arbitrary dimension.

\subsection{Poset of proper subspaces} 
Let $\tT(V)$ denote the poset of proper subspaces of $V$, and $\tT^\circ(V)=T(V)-\set{0}$. The geometric realization of the latter poset is usually called the Tits building associated to the general linear group $\GL(V)$. Its homotopy type  is described by the   Solomon--Tits theorem.
\Pro{\cite{S69} \label{ST} Assume that $\dim V\geq 1$. Then $\tT^\circ(V)$ is $(\dim V -2)$--spherical.}
\noindent In his work on discrete series representations of the classical groups Lusztig studied the poset $\cC_V\tT(V)$. In \cite[Theorem 1.9]{L74} he proves a homological version of the following. 
\Pro{\label{affine} The coset poset  $\cC_V\tT(V)$ is  $(\dim V-1)$--spherical.} 
\noindent This result is a special case of \cite[Proposition 11]{B00} where it is extended to cosets of proper subgroups of a  solvable group.

Let $W$ be a subspace of $V$ where $\codim\; W=1$. Recall that we have the following maps from Proposition \ref{H-eq} and \ref{thetaH}:
$$
i_W: \cC_W\tT(V) \rightarrow  \cC_W \tT(W) \;\; \text{ and } \;\;
\theta_W: \cC_V \tT(V)^{\vee W} \rightarrow \cC_W \tT(V)
$$
where we identified $\tT(V)_{\cap W}$ with $\tT(W)$ in the first map.
Fix $v\notin W$ and  define a map of posets
\begin{equation}\label{theta_v}
\theta_v :\tT(V)^{\vee W}\rightarrow \cC_W\tT(W)
\end{equation}
 where $\theta_v(A)=(-v+A)\cap W$. This map is well-defined since it can be factored as  
\begin{equation}\label{thetav}
 \begin{tikzcd}
\tT(V) ^{\vee W} \arrow{r}{\theta_v} \arrow{d}{-v} & \cC_W\tT(W)\\
-v+\tT(V) ^{\vee W}\arrow{r}{\theta_W} & \cC_W\tT(V) \arrow{u}{i_W} \\
\end{tikzcd}
\end{equation}
where   $\theta_{W}$ is  restricted on the subposet
$
-v+\tT(V) ^{\vee W} \subset  \cC_V \tT(V)^{\vee W}
$
of cosets of the form $-v+A$.
 
\Pro{\label{rel}Assume that $\dim V\geq 1$.  The map $\theta_v :\tT(V)^{\vee W}\rightarrow \cC_W\tT(W)$ is a weak equivalence. In particular $\tT(V)^{\vee W}$ is $(\dim V-2)$--spherical.} 
\Proof{ Its homotopy inverse is given by the map
$$
s_v:\cC_W\tT(W)\rightarrow \tT(V)^{\vee W}
$$ 
defined by $s_v(w+A)=\Span{v+w,A}$. One checks that $\theta_vs_v(w+A) \supset w+A$, and $s_v\theta_v(B)\subset B$. 
Therefore $s_v$ is the homotopy inverse by the basic result   described in Section \ref{homotopy_poset} and the homotopy type of the coset poset is given in Proposition \ref{affine}. 
}

There is a variation of this result when $\dim V\geq 2$ whose proof is essentially the same. 
Let $U$ be a subspace of $W$ with $\dim U=1$. 
Let $\tT(V)_{\wedge U}$ denote the poset $\set{A\in \tT(V)|\; A\cap U=0}$.
 We will denote the intersection $\tT(V)^{\vee W}\cap \tT(V)_{\wedge U}$ simply by $\tT(V)^W_U$.   
Consider the restriction of $\theta_v$ to the subposet $\tT(V)^W_U$. The image of the restricted map $\theta'_v$ lies  in $\cC_W\tT(W)_{\wedge U}$.  Therefore we  have a diagram 
\begin{equation*}
\begin{tikzcd}
\tT(V)^{\vee W} \arrow{r}{\theta_v} & \cC_W\tT(W) \\
\tT(V)^W_U \arrow[u,hook] \arrow{r}{\theta'_v} & \cC_W\tT(W)_{\wedge U}\arrow[u,hook]
\end{tikzcd}
\end{equation*}
where the vertical maps are natural inclusions.
 The map $\theta'_v$ is also a weak equivalence whose  homotopy inverse is given by the restriction of $s_v$ to $\cC_W\tT(W)_{\wedge U} $. We record this result.

\Pro{\label{rel2}Assume that $\dim V\geq 2$.  The map $\theta_v' :\tT(V)^{W}_U\rightarrow \cC_W\tT(W)_{\wedge U}$ is a weak equivalence.} 
 
When $\dim V =2$ the poset $\tT(V)^{W}_U$ consists of one-dimensional subspaces different than $W$. There are $p$ many of such subspaces. The coset poset $\cC_W\tT(W)_{\wedge U}$ consists of the cosets of the zero subspace in $W$. Again there are $p$ many of such cosets.  
The poset of subspaces $\tT(V)^{W}_{U}$ will appear as fibers of certain maps between coset posets. Next we show that this poset is spherical.

\Thm{\label{general}Assume that $\dim V\geq 2$. Let $U\subset W$ be subspaces in $V$ and $\codim\, W=\dim\, U=1$. Then $ \tT(V)^W_U$ is $(\dim V-2)$--spherical. }  

\begin{proof} 
By Proposition \ref{rel2} the poset $\tT(V)^W_U$ is weakly equivalent to the coset poset $\cC_W\tT(W)_{\wedge U}$.
The statement of the theorem holds when $\dim V=2$ since the resulting space is a disjoint union of points. 
For larger values of $\dim\,V$  we will proceed by induction on the dimension. 
Assume that the statement of the theorem holds for vector spaces of dimension less than $\dim V$.   We will decompose $\cC_W\tT(W)_{\wedge U}$ with respect to a subspace $L$  of $W$ of codimension one which contains $U$.  
By Corollary \ref{exact-n} it suffices to show that $\cC_L\tT(W)_{\wedge U}$ and $\cC_W\tT(W)_U^L$  are $(\dim W-2)$--spherical, and construct the map $s$. We start with the homotopy type of the first one. Since the poset $\tT(W)_{\wedge U}$ is $L$--stable Proposition \ref{H-eq} implies that 
$$i_L:\cC_L\tT(W)_{\wedge U}\rightarrow \cC_L(\tT(W)_{\wedge U})_{\cap L}=\cC_L\tT(L)_{\wedge U}$$
is a weak equivalence and by Proposition \ref{rel2} there is a weak equivalence $\cC_L\tT(L)_{\wedge U} \simeq \tT(W)^L_U$
   which is $(\dim W-2)$--spherical by the induction hypothesis. 
Next we determine the homotopy type of $\cC_W\tT(W)_U^L$ by applying Theorem \ref{Quillen} to the natural map $\epsilon:\cC_W\tT(W)_U^L\rightarrow \tT(W)_U^L$ defined by $w+A\mapsto A$. We need to consider the posets
$
(\tT(W)_U^L)_{<A} 
$ and
$\overf{\epsilon}{A}
$
and show that both are spherical. The first one can be identified as follows:
$$
(\tT(W)_U^L)_{<A}=(\tT(W)^{\vee L})_{<A}=\tT(A)^{\vee L\cap A}
$$
where the latter  is spherical by Proposition \ref{rel}. The fiber can be identified as
$$
\overf{\epsilon}{A}=\cC_W(\tT(W)_U^L)_{\geq A}=\cC_W(\tT(W)_{\wedge U})_{\geq A} 
$$
and there is a map of posets
$$
\alpha:\cC_W(\tT(W)_{\wedge U})_{\geq A}\rightarrow \cC_{W/A}\tT(W/A)_{\wedge UA/A}
$$
defined by $w+B\mapsto w+BA/A$. We claim this is a weak equivalence. To see this let $\bar B$ denote a coset in the image and consider the fiber $\underf{\alpha}{\bar B}$. The fiber is contractible since it has a terminal object. Hence  Theorem \ref{QuillenThA} implies that $\alpha$ is a weak equivalence. Finally, again using Proposition \ref{rel2} and the induction hypothesis we see that the fiber is spherical.
 Therefore both spaces satisfy the requirements in Theorem \ref{Quillen} and we conclude that $\cC_W\tT(W)_U^L$ is $(\dim W-2)$--spherical.
It remains to define the map $s:\cC_L \tT(W)_{\wedge U}  \rightarrow \tT(W)^L_U$. According to Theorem \ref{exact} this map needs to satisfy the following two properties: $\theta_Ls \simeq \ast$ and the composite
$$
\theta_L(-v)s:\cC_L \tT(W)_{\wedge U}  \stackrel{s}{\longrightarrow} \tT(W)^L_U \stackrel{-v}{\longrightarrow} -v+ \tT(W)^L_U \stackrel{\theta_L}{\longrightarrow}   \cC_L\tT(W)_{\wedge U}
$$
is homotopic to the identity map for some fixed element $v\in W- L$. Here we restricted the map $\theta_L:\cC_{W} \tT(W)_U^L  \rightarrow\cC_{L}\tT(W)_{\wedge U} $ to the subposet of cosets of the form $-v+A$ where $A\in \tT(W)^L_U$.
Our candidate for $s$ is closely related to the homotopy inverse $s_v$ of the map   
$$
\theta'_v:\tT(W)^L_U\rightarrow \cC_L\tT(L)_{\wedge U}
$$ 
considered in Proposition \ref{rel2}.  Note that by diagram \ref{thetav} we have $\theta_v'=i_L\theta_L(-v)$. Let $j$ be a homotopy inverse to $i_L$. We claim $s= s_vi_L$ is the map we need.
We have $\theta_L(-v)s \simeq j\theta_v's_vi_L \simeq ji_L \simeq \idy$. It remains to check $\theta_Ls = \theta_L s_v i_L \simeq \ast$. For this recall that  $s_v$ sends a coset $l+A$ to the subspace $\Span{v+l,A}$.
 By direct verification we see that $\theta_L s_v i_L(l+A)$ is the trivial coset $\Span{v+l,A\cap L}$ which contains $0$, yielding a contracting homotopy $\theta_Ls\simeq \ast$.  
\end{proof}

\subsection{Poset of isotropic subspaces}  Let $\bi$ be an alternating bilinear form on $V$. This means $\bi(v,v')=-\bi(v',v)$ for all $v,v'$ in $V$.   The orthogonal complement of a  subspace $A\subset V$   is defined by $A^\perp=\set{w\in V|\;\bi(w,a)=0\;\forall\;a\in A}$. For an element $v$ we write $v^\perp$ to denote the orthogonal complement of the subspace $\Span{v}$. We assume $\bi$ non-degenerate i.e. $V^\perp=0$. A subspace is called isotropic if $A\subset A^\perp$. 
Let  $\iI(V)$ denote the collection of isotropic subspaces.
Let $H(V)$ denote the Heisenberg group associated to the vector space $V$ with the bilinear form $\bi$ as introduced in section \ref{Heisenberg}. 
We can regard $\iI(V)$ as a collection of abelian subgroups for $H(V)$ under the inclusion $I\rightarrow H(V)$ defined by $a\mapsto (a,0)$.
Then the centralizer of $I$ in $H(V)$ can be identified with the subgroup $H(I^\perp)$.

\Thm{\label{main1} There is a weak equivalence
$$
\cC_{H(V)} \iI(V) \simeq \bigvee^{d(p,r)} \sS^r
$$
 where $2r=\dim V$ and $d(p,r)$ is defined in Corollary \ref{euler}.
}

The proof of this theorem will occupy the rest of this section. 
We fix a symplectic basis $b=\set{x_1,x_2,\cdots,x_r,\bar{x}_1,\bar{x}_2,\cdots,\bar{x}_r}$ of $V$ with  $\bi(x_i,\bar{x}_j)=1$ when $|i-j|=0$ and zero otherwise.
Let us fix $x=x_r$ and $\bar{x}=\bar{x}_r$ for the rest. 
We will apply Corollary \ref{exact-n} with $H=H(x^\perp)$ which is a normal subgroup of $H(V)$ with a quotient  isomorphic to a cyclic group  of order $p$. The coset space $\cC_{H(V)}\iI(V)$ is the homotopy colimit of the diagram
\begin{equation*} 
\begin{tikzcd}
\cC_{gH } \iI(V) & \cdots &  \cC_{g^pH } \iI(V) \\
 &\arrow{ul}{\theta_{gH}} \cC_{H(V)}\iI(V)^{\vee H} \arrow[ur,"\theta_{g^pH}"']&
\end{tikzcd}
\end{equation*}
where $g=(\bar x)$ and $H=H(x^\perp)$. Note that $g^p=(p\bar x)=(0)$ is the identity element of $H(V)$.
We have two  main objectives: 1) construct the map $s$, and 2) show that $\cC_{H(V)} \iI(V)^{\vee H}$ is $(r-1)$--spherical.

\subsection{Construction of $s$}  We will construct a map of posets  
$$s: \cC_H\iI(V)\rightarrow \iI(V)^{\vee H} $$
such that $\theta_Hg^{-1}s\simeq \idy$ and $\theta_Hs\simeq \ast$ as required by Corollary \ref{exact-n}.
The quotient map $p:x^\perp\rightarrow x^\perp/\Span{x}$ induces a weak equivalence
\begin{equation}\label{p}
p:\cC_{H } \iI(x^\perp) \rightarrow \cC_{H(x^\perp/\Span{x})}\iI(x^\perp/\Span{x})
\end{equation}
since   the fibers canonically contract to a terminal object. Here $H(x^\perp)/\Span{x}$ is identified with $H(x^\perp/\Span{x})$. Moreover by Proposition \ref{H-eq} the map
\begin{equation}\label{iH}
i_{H}:\cC_{H }\iI(V)\rightarrow \cC_{H }\iI(V)_{\cap H }=\cC_{H }\iI(x^\perp)
\end{equation}
is a weak equivalence since $\iI(V)$ is $H$--stable. Let $\phi$ denote the composition 
\begin{equation}\label{phi}
\phi=pi_{H}: \cC_{H }\iI(V)\rightarrow \cC_{H(x^\perp/\Span{x})}\iI(x^\perp/\Span{x})
\end{equation} 
 and let $\bar\theta_H$ denote the composition $\phi\theta_H$ where $\theta_{H }:\iI(V)^{\vee H }\rightarrow \cC_{H }\iI(V)$.

\Lem{\label{sbar}
Suppose that the map of posets
\begin{equation}\label{bartheta}
\bar\theta=\bar\theta_{H }g^{-1}:\iI(V)^{\vee H(x^\perp)}\rightarrow \cC_{H(x^\perp/\Span{x})}\iI(x^\perp/\Span{x}) 
\end{equation}
has an inverse denoted by $\bar{s}$ such that $\bar{\theta}_{H } \bar s\simeq \ast$. Then $s=\bar s \phi$ satisfies the required properties in Corollary \ref{exact-n}.
}
\Proof{
If $\varphi$ is a homotopy inverse of $\phi$ then multiplying $\bar\theta \bar s=\idy$ by $\varphi$ on the left and by $\phi$ on the right we get 
$$
\idy\simeq \varphi \phi = \varphi \bar\theta\bar s \phi =\varphi  \phi\theta_{H }g^{-1}  s\simeq \theta_{H }g^{-1}  s.
$$
 Moreover $\theta_{H}s\simeq \ast$ since  $\phi\theta_{H}s\phi\simeq\bar\theta_{H}\bar s \simeq \ast$.
}
We proceed with  the construction of $\bar s$. 
The quotient $x^\perp/\Span{x}$ can be identified with the subspace $Z$  spanned by $b-\set{x,\bar{x}}$. Let $j:x^\perp/\Span{x}\rightarrow Z\subset x^\perp$ denote this identification. There is a commutative diagram
\begin{equation*}
\begin{tikzcd}
x^\perp \arrow{r}{k} \arrow{d}{p} & x^\perp\\
x^\perp/\Span{x} \arrow{ur}{j} & 
\end{tikzcd}
\end{equation*}
where  $k(u)=u-\bi(u,\bar{x})x$. 
There is a corresponding  commutative diagram 
\begin{equation}\label{dmapH}
\begin{tikzcd}
H(x^\perp) \arrow{r}{k} \arrow{d}{p} & H(x^\perp)\\
H(x^\perp/\Span{x}) \arrow{ur}{j} & 
\end{tikzcd}
\end{equation}
of Heisenberg groups, where the  maps are still denoted by the same letters. 
We define a map of posets
\begin{equation}\label{inverseH}
\bar s:  \cC_{H(x^\perp/\Span{x})}\iI(x^\perp/\Span{x}) \rightarrow \iI(V)^{\vee H}
\end{equation}
which sends a coset $(w,t)A\subset j(H(x^\perp/\Span{x}))$ to the subgroup 
$$
\Span{ (\bar{x}+tx+w),(\bi(w,a)x+a)  |\;a\in A  }.
$$ 
Let us check that this definition gives an abelian subgroup of $H$.
Note that the commutator   of two elements in $H(V)$ is given by 
$$[(v_1),(v_2)]=(v_1,0)(v_2,0)(-v_1,0)(-v_2,0)=(v_1+v_2,\bi(v_1,v_2))(-v_1-v_2,\bi(v_1,v_2))= (0,2\bi(v_1,v_2)).$$
Therefore $\bar{s}$ applied to a coset gives an abelian subgroup since 
$$\bi(\bar{x}+tx+w,\bi(w,a)x+a)=\bi(w,a)\bi(\bar{x},x)+\bi(w,a)=0.$$ 
Next we prove a key lemma.
 
\Lem{\label{invpH}The map $\bar s$ is the inverse of $\bar \theta$ as defined in \ref{bartheta}. It satisfies $\bar{\theta}_{H} \bar s\simeq \ast$.}
\Proof{First we show that $\bar s \bar\theta(A)=A$ for all $A$ in $\iI(V)^H$. Such isotropic subspaces satisfy
 $A+x^\perp=V$ where we regard $A$ as a subgroup of $H(V)$ via the   mapping $a\mapsto (a,0)$. We have $\bar{x}+u\in A$ for some $u\in x^\perp$ since $A$ and $x^\perp$ span the whole space $V$.
Then  $A=\Span{\bar{x}+u,A\cap x^\perp}$. Recall that $\bar\theta=\bar{\theta}_Hg^{-1}$ sends a space $A$ to the quotient $((g^{-1}A)\cap H)/\Span{x}$ where $H=H(x^\perp)$ and $g=(\bar x)$. Although the function $\bar\theta$ is defined on subspaces it will be convenient to think that it is defined as follows on  vectors $v$ such that $\Span{v} +x^\perp =V$: 
$$\bar{\theta}(v)=p((\bar{x})^{-1}(v)\cap H )=(p(v-\bar{x}),\bi(v,\bar{x}))  
$$
 where $p:x^\perp\rightarrow x^\perp/\Span{x}$. Then we can write
$$j\bar{\theta}(v)=(k(v-\bar x ),\bi(v,\bar{x})) $$
 by the commutativity of \ref{dmapH}.
In particular, we have $j\bar\theta (\bar{x}+u)= (k(u),\bi(u,\bar{x}))$. More generally for $a\in A\cap x^\perp$ we compute
\begin{align*}
j\bar\theta(\bar{x}+u+a)&= (k(u+a),\bi(u+a,\bar{x}))\\
&= (k(u),\bi(u,\bar{x}))(k(a),\bi(a,\bar{x})+\bi(a,u))\\
&= (k(u),\bi(u,\bar{x}))(k(a),\bi(a,\bar{x}+u))\\
&= (k(u),\bi(u,\bar{x}))(k(a),0)
\end{align*}
where in the last step we have used the fact that $a$ and $\bar{x}+u$ are elements of an isotropic space $A$ i.e. $\bi(a,\bar{x}+u)=0$. Therefore the composition $j\bar{\theta}$ maps $A$ to the coset
$$
A=\Span{\bar{x}+u,A\cap x^\perp}\;\; \mapsto\;\; (k(u),\bi(u,\bar{x}))k(A)
$$ 
in $j(H(x^\perp/\Span{x}))$.
Applying $\bar s$ to this coset gives a subgroup generated by 
$(\bar{x}+\bi(u,\bar{x})x+k(u))=(\bar{x}+u)$ and 
$$(\bi(k(u),k(a))x+k(a))=(\bi(u,a)x+k(a))=(\bi(\bar x+u,a)x+a)=(a)$$
where we have used the fact that $k$ is defined by $k(u)=u-\bi(u,\bar{x})x$ and respects the bilinear form.
Therefore
$$
\bar s \bar\theta(A)=\Span{(\bar{x}+u),(a)|\;a\in A\cap x^\perp}
$$
which is exactly equal to $A$ when regarded as a subgroup of $H(V)$.

For the converse let $(w,t)A$ be a coset in the image of $j$.  
After applying $\bar{s}$  consider the element $(\bar{x}+tx+w+\bi(w,a)x+a,\bi(\bar x,a))$ obtained as the product of the two generators $(\bar{x}+tx+w)$ and $(\bi(w,a)x+a)$.
To see what $\bar{\theta}$ does to $\bar{s}((w,t)A)$ it suffices to check its effect on this element  
\begin{align*}
(-\bar{x})(\bar{x}+tx+w+\bi(w,a)x+a,\bi(\bar x,a))&=(tx+w+\bi(w,a)x+a,t+\bi(w,a))\\
&\equiv (w+a,t+\bi(w,a)) \mod (x)\\
&=  (w,t)(a,0) \mod (x)
\end{align*}
which implies that $\bar\theta\bar s$ is the identity. 

Finally, the last statement follows from $\bar\theta_{H }\bar s((w,t)A)\supset (0)$ since in effect this composite is computed by intersecting the subgroup $\bar s((w,t)A)$  with $H(x^\perp)$ and taking the quotient by $\Span{x}$. Therefore this composition  contracts to the constant map at the identity element.
}

\subsection{Proof of Theorem \ref{main1}}
Restricting the domains of the maps $i=i_{H}$ in \ref{iH} and $p$ in \ref{p}, we obtain the following maps 
\begin{equation}\label{ip}
\cC_{H}\iI(V)^{\vee H} \stackrel{i}{\rightarrow}\cC_{H}\iI(x^\perp)_{\wedge\Span{x}}  \stackrel{p}{\rightarrow}
\cC_{H(x^\perp/\Span{x})}\iI(x^\perp/\Span{x})
\end{equation}
where $\iI(x^\perp)_{\wedge\Span{x}}$ denotes the poset of subspaces $A$ in $\iI(x^\perp)$ such that $A\cap \Span{x}=0$.
Let us consider the fibers of each of these maps. Since both maps are equivariant with respect to the left action of $H$ it suffices to consider the fiber over  a coset of the form $B$. Other cosets $\alpha B$ are isomorphic to this one under the action of $\alpha \in H$.

\Lem{\label{fib1} Let $B\in \iI(x^\perp/\Span{x})$. The fiber  $\underf{p}{B} $ is $(\dim B)$--spherical.}
\Proof{Lifting $B$ to an isotropic subspace $\tilde{B}=\Span{j(B),x}$ in $x^\perp$, we have
$$ \underf{p}{B} = \cC_{\tilde{B}}\tT(\tilde{B})_{\wedge \Span{x}}\simeq \tT(\tilde{B}\oplus \FF_p)^{ \tilde{B}}_{\Span{x}}$$
where the weak equivalence follows from Theorem \ref{general} with $V=\tilde{B}\oplus \FF_p$ and hence it is $(\dim(\tilde{B}\oplus \FF_p)-2)$--spherical as desired. 
}

\Lem{\label{fib2}Let $B\in \iI(x^\perp)_{\wedge\Span{x}}$.
The fiber  $\overf{ i}{B} $ is $(r-\dim B-1)$--spherical.
} 
\Proof{ Let $B^\perp$ denote the orthogonal complement of $B$ in $V$ with respect to $\bi$. Note that $B^\perp$ contains $x$ since $B\subset x^\perp$.
We have $$\overf{ i}{B}=(\iI(V)^{\vee H})_{\geq B}\cong \iI(B^\perp/B)^{\vee H(B^\perp\cap x^\perp)/B}$$
where the last poset is isomorphic to $\cC_{H(B^\perp\cap x^\perp/\Span{B,x}) }\iI(B^\perp\cap x^\perp/\Span{B,x})$ by Lemma \ref{invpH}. The dimension of a maximal isotropic subspace in  $B^\perp\cap x^\perp/\Span{B,x}$ is $r-\dim B-1$. }
 
Now we can finish the proof of Theorem \ref{main1} by collecting the results obtained so far. 
We will do so by induction on $r$ where $2r=\dim V$. When $r=2$ the coset poset $\cC_{H(V)}\iI(V)$ is a one-dimensional connected space hence it is $1$--spherical.
We already constructed the section $s$ in Lemma \ref{sbar} and Lemma \ref{invpH}.
Using the weak equivalence $\phi$ in \ref{phi}   and induction we see that $\cC_{H } \iI(V)$ is $(r-1)$--spherical. 
It remains to prove that
$$\xX=\cC_{H(V)}\iI(V)^{\vee H }$$ is $(r-1)$--spherical. We will use Lemma \ref{fib1} and Lemma \ref{fib2}.
First consider the map 
$$
p:\xX'\rightarrow \yY'
$$
in \ref{ip}. By the induction assumption $\yY'$ is $(r-1)$--spherical.  Let $B\in \iI(x^\perp/\Span{x})$. The fiber $\underf{p}{B}$ is $(\dim B)$--spherical by Lemma \ref{fib1}, and $\yY'_{>B}=\cC_{H(B^\perp/B)}\iI(B^\perp/B)$ is $(r-1-\dim B)$--spherical by inspection. Then $\xX'$ is $(r-1)$--spherical by Theorem \ref{Quillen} applied to the opposite of the map $\xX'\rightarrow \yY'$ (original statement in \cite[Theorem 9.1]{Q78}). We turn to the other map in \ref{ip}, and denote it simply by
 $$i:\xX\rightarrow \xX'.$$
  Let $B\in \iI(x^\perp)_{\wedge\Span{x}}$.
 In this case $\xX'_{<B}=\cC_B\tT(B)_{\wedge \Span{x}}$  is $(\dim B-1)$--spherical by Proposition \ref{rel} and $\overf{i}{B}$ is $(r-\dim B-1)$--spherical by Lemma \ref{fib2}. Therefore $\xX$ is $(r-1)$--spherical by Theorem \ref{Quillen} and decomposes as a wedge of    spheres where the number of spheres is equal to $d(p,r)$ by Corollary \ref{euler}.
This concludes the proof of Theorem \ref{main1}.

\section{Classifying space for commutativity} \label{secapp}

In \cite{ACT12} a natural filtration $\set{B(q,G)}_{q\geq 2}$ of the classifying space $BG$ is introduced. These spaces can be described as homotopy colimits of classifying spaces and coset posets naturally occur in the study of such objects. In this section we apply our main result Theorem \ref{main0} to $B(2,G)$ when $G$ is an extraspecial $p$--group.

\subsection{The space $B(2,G)$} As a simplicial set  $B(2,G)$ has $n$--simplices  given by the set of group homomorphisms $\Hom(\ZZ^n,G)$ and the simplicial structure is induced from $BG$.  
Let $E(2,G)\rightarrow B(2,G)$ denote the pull-back of the universal principal $G$--bundle $EG\rightarrow BG$ along the natural inclusion $B(2,G)\subset BG$. 
It is known that the natural map
$$
\hocolim{\aA(G)} B \rightarrow B(2,G)
$$
is a weak equivalence \cite[\S 4]{ACT12}. As a consequence of Theorem \ref{main0} we determine the universal cover of $B(2,G)$ when $G$ is an extraspecial $p$--group.

\Thm{\label{main2}
Let $E$ be an extraspecial $p$--group. There is a fibration sequence
$$
\bigvee^{d(p,r)} \sS^r \rightarrow B(2,G) \rightarrow B\pi
$$
where $2r\geq 4$ is the rank of the Frattini quotient of $E$, the number $d(p,r)$ is defined in Corollary \ref{euler}, and $\pi$ is the kernel of the multiplication map
$$
1\rightarrow \pi \rightarrow E\times E \rightarrow E/[E,E]\rightarrow 1.
$$
}
\Proof{
By Proposition \ref{colim_ext} the fundamental group of $B(2,G)$ is isomorphic to $\pi$. The universal cover of $B(2,G)$ is identified with the coset poset $\cC_\pi\aA(G)$ and hence with $\cC_{HV}\iI(V)$ by Proposition \ref{diag_covers}.
}

Let us introduce a variant of $B(2,G)$. Let $V$ be a vector space over $\FF_p$ with a non-degenerate alternating bilinear form $\bi$. We define a simplicial set $B(\bi,V)$ whose set of $n$--simplices is the set of $n$--tuples $(v_1,v_2,\cdots,v_n)$ where $\bi(v_i,v_j)=0$ for all $0\leq i,j\leq n$. The simplicial structure is induced via the inclusion $B(\bi,V)\subset BV$. Similar to $B(2,G)$ we can describe this space as a homotopy colimit. There is a natural weak equivalence
$$
\hocolim{\iI(V)} B \rightarrow B(\bi,V)
$$
induced by the inclusions $BI\subset B(\bi,V)$ where $I$ is an isotropic subspace.
Let $E(\bi,V)\rightarrow B(\bi,V)$ denote the pull-back of the universal bundle $EV\rightarrow BV$ along the inclusion $B(\bi,V)\subset BV$. There is a weak equivalence
$$
E(\bi,V) \simeq \cC_{V} \iI(V)
$$
since both spaces are the homotopy fibers of the map $B(\bi,V)\rightarrow BV$.

\begin{pro}The projection $E\rightarrow V$ induces a diagram
$$
\begin{tikzcd}
E\ZZ/p \arrow{r} \arrow{d} & E(2,E) \arrow{r} \arrow{d} & E(\bi,V) \arrow{d} \\
B\ZZ/p \arrow{r} & B(2,E) \arrow{r} & B(\bi,V)
\end{tikzcd}
$$
where the horizontal maps are fibration sequences. Moreover, there is a fibration sequence
$$
\bigvee^{d(p,r)} \sS^r \rightarrow B(\bi,V) \rightarrow B(H(V))
$$
where $H(V)$ is the Heisenberg group.
\end{pro} 
\Proof{
It is a direct verification to show that $E(2,E)\rightarrow E(\bi,V)$ is a Kan fibration by checking the horn lifting property \cite[pg. 11]{GJ99}. The idea is that it is always possible to lift an isotropic subspace $I$ to an  abelian subgroup of $E$. We think of the zero element    as the base point $\ast$ of $E(\bi,V)$ and over this point the fiber is $E\ZZ/p$ by definition of the map. A similar argument works for the map $B(2,E)\rightarrow B(\bi,V)$. The last statement follows from Theorem \ref{main2} and  Corollary \ref{covering} where we showed that the fundamental group of the homotopy colimit of $BI$ where $I\in \iI(V)$ is the Heisenberg group $H(V)$.
}

\subsection{Subdivision} Let $X$ be a simplicial set. We denote its subdivision by $\sd\,X$. It comes with a map $l:X\rightarrow \sd\,X$ called the last vertex map. This map is a weak equivalence \cite[pg. 193]{GJ99}. For certain simplicial sets there is a nicer description of the subdivision construction. Let $\Ba(X)$ denote the poset of non-degenerate simplices of $X$ ordered under the face relation: $x\leq y$ if $x$ is a face of $y$. If any non-degenerate $n$--simplex of $X$ has $n+1$ distinct vertices then $\sd\,X$ is isomorphic to the nerve of the poset $\Ba(X)$, called the Barratt nerve \cite[Lemma 2.2.11]{WJR13}. Some examples of such simplicial sets are $EG$, $E(2,G)$, $E(\bi,V)$,...
We identify the subdivision of such a simplicial set with its Barratt nerve.
There is a map of posets
\begin{equation}\label{subtheta}
\theta: \sd\, E(\bi,V) \rightarrow \cC_V \iI(V)
\end{equation}
defined by sending a non-degenerate simplex $(v,a_1,a_2,\cdots,a_n)$ to the coset $v+\Span{a_1,a_2,\cdots,a_n}$.

\subsection{Transitionally commutative bundles} A principal $G$--bundle $p:E\rightarrow X$ is called transitionally commutative if there exists an open cover $\set{U_j}_{j\in J}$ of $X$ such that the restricted bundle $p|_{U_j}$ is trivial for all $j\in J$ and the transition functions commute when simultaneously defined. Two transitionally commutative bundles $p_0$ and $p_1$ are said to be isomorphic if there exists a transitionally commutative principal $G$--bundle $p:E\rightarrow X\times[0,1]$ such that $p|_{X\times 0}=p_0$ and $p|_{X\times 1}=p_1$.

\Thm{\cite[Theorem 2.2]{AG15} Let $X$ be a finite CW-complex and $G$ a  Lie group. There is a one-to-one correspondence
between the set of homotopy classes of maps $[X,B(2,G)]$ and the set of isomorphism classes of transitionally commutative principal $G$--bundles.
}

The bundle $E(2,G)\rightarrow B(2,G)$  is the universal example of a transitionally commutative $G$--bundle in the sense that any transitionally commutative $G$--bundle over $X$ is isomorphic to a pull-back bundle $f^*(E(2,G))\rightarrow X$ for some map $f:X\rightarrow B(2,G)$.

\Cor{Let $E$ be an extraspecial $p$--group.
Isomorphism classes of transitionally commutative principal $E$--bundles over $\sS^r$ are given by
$$
[\sS^r,B(2,E)]=[\sS^r,\bigvee^d \sS^r] =\ZZ^d
$$
where $2r\geq 4$ is the Frattini quotient of $E$.
}   
\Proof{
Since $r\geq 2$ any map $\sS^r\rightarrow B(2,E)$ lifts to the universal cover whose homotopy type is a wedge of spheres as determined by Theorem \ref{main2}. The result follows from the Hurewicz theorem.
}

This result says that although any principal $E$--bundle over $\sS^r$ is trivial for $r\geq 2$ such bundles cannot be trivialized through transitionally commutative principal $E$--bundles.

\subsection{A non-trivial bundle} Next we construct a non-trivial transitionally commutative principal $E$--bundle over the $r$--sphere. 
Let $\set{x_i,\bar x_i}_{i=1}^r$ denote a symplectic basis for $V$. Set $x= \sum_{i=1}^r (x_i+\bar x_i)$ and $\bar x = 0$. Consider the join of two point discrete posets
$$
J= \set{x,\bar{x}} \ast \set{x_1,\bar{x}_1} \ast \cdots \set{x_r,\bar{x}_r}
$$
whose geometric realization is homeomorphic to the $r$--sphere. 
For each object $j$ of $J$ let us fix a lift $\tilde j$ with respect to the projection $E\rightarrow V$.
Identifying $J$ with its nerve  we define a map of simplicial sets  
$$
\tau:J\rightarrow B(2,E)
$$
by sending an $n$--simplex $j_0\rightarrow j_1 \rightarrow \cdots \rightarrow j_n$ to the $n$--simplex 
$$((\tilde j_0)^{-1}\tilde j_1 , (\tilde j_1)^{-1}\tilde j_2,\cdots,  (\tilde j_{n-1})^{-1}\tilde j_{n}).$$ 
Our goal is to show that the class represented by $\tau$ in the homotopy group $\pi_rB(2,E)$ is non-trivial.  
For this we introduce another map which will turn out to be closely related to $\tau$.
Let us define a map of posets
$$
\tilde\tau: \sd\, J \rightarrow \cC_{H(V)} \iI(V) 
$$
by sending a chain $j_1<j_2<\cdots<j_s$ of objects in the poset $J$ to the coset $(j_1)\Span{\,j_2-j_1, j_3-j_1,\cdots, j_s-j_1\,}$ in $H(V)$.

\Lem{\label{nontrivial} The class represented by $\tilde\tau$ in the homotopy group $\pi_r\,  \cC_{H(V)} \iI(V)$ is non-trivial.} 
\begin{proof}
We first show that $\tilde\tau$ factors as
$$
\sd\, J \stackrel{\bar \tau}{\longrightarrow} \iI(V')^{\vee H(x_{0}^\perp)}  \stackrel{\bar \theta}{\longrightarrow}    \cC_{H(V)} \iI(V) 
$$ 
where $V'=\Span{x_{0},\bar x_{0}}\oplus V$ is a symplectic vector space of dimension $2(r+1)$, $x_{0}^\perp$ is the orthogonal complement in $V'$ with respect to the standard symplectic form, and $\bar \theta$ is the isomorphism in \ref{bartheta} defined with respect to the vector $x_0$ as follows
$$
A \mapsto (\,((\bar x_{0})^{-1} A) \cap H(x_{0}^\perp) \,)/\Span{x_{0}}.
$$ 
To define $\bar\tau$   we start with the canonical chamber in the Tits building associated to $\Sp(V')$ defined by a map of posets
$$
c: \sd\,J'\rightarrow \iI(V')^\circ
$$
where $J'=\set{x_0,\bar{x}_0}\ast \set{x_1,\bar x_1}\ast \cdots \ast \set{x_{r},\bar x_{r}}$ and a chain of objects  $i_1<i_2<\cdots <i_t$ maps to the subspace $\Span{i_1,i_2,\cdots,i_t}$. This map is non-trivial in homology by the Solomon-Tits theorem and hence gives a non-trivial map in homotopy. Let $g$ denote the symplectic transformation on $V'$ defined by 
$$
x_{0} \mapsto \sum_{j=0}^{r} (x_j+\bar x_j) ,\; \bar x_{0}\mapsto \bar x_{0},\; x_i \mapsto x_i+\bar x_{0},\; \bar x_i \mapsto \bar x_i+\bar x_{0},\;\;\; 1\leq i\leq r.
$$ 
Let $\hat g$ denote the map of posets $\iI(V')^\circ\rightarrow \iI(V')^\circ$ defined by $A\mapsto g(A)$. 
The image of the composite $\hat gc$ lies in the poset $\iI(V')^{\vee H(x_{0}^\perp)}$.
We define $\bar{\tau}=\hat g c h$ where  $h$ is the isomorphism between the subdivisions $\sd\,J\rightarrow \sd\,J'$ induced by the map of posets $J\rightarrow J'$ defined by $x\mapsto  x_{0}$, $\bar x\mapsto\bar x_{0}$ and the remaining $x_i$, $\bar{x}_i$  for $1\leq i\leq r$ are fixed. Note that by construction $\bar \tau$ is not homotopic to the constant map. It remains to check that $\tilde\tau=\bar\theta \bar \tau$. Let $j_1<j_2<\cdots<j_s$ be an object of $\sd\, J$. Under  $\bar \tau$ it maps to a subgroup $\Span{j_1',j_2',\cdots,j_s'}$ where $j_1'=j_1+x_0+\bar x_0$ if $j_1=x$ or $j_1'=j_1+ \bar x_0$ if $j_1\not= x$, and the remaining generators are given by $j_k'=j_k+\bar x_0$ for $1<k\leq s$.  Now the elements in the image of $\bar \theta$ come from elements of the form $\alpha_1j_1'+\alpha_2j_2'+\cdots +\alpha_sj_s'$ where $\sum_k \alpha_k \equiv 1$ mod $p$. We can write
\begin{align*}
\sum_{k=1}^s \alpha_kj_k' &= j_1+ \sum_{k=2}^s \alpha_k (j_k'-j_1) \\
&\equiv \bar x_0 + j_1+ \sum_{k=2}^s \alpha_k (j_k-j_1)  \mod \Span{x_0}
\end{align*}
and thus the image under $\bar{\theta}$ is the coset $(j_1)\Span{j_2-j_1,\cdots,j_s-j_1}$ in $H(V)$.
\end{proof}

\Thm{\label{nontrivial_construction} Assume $r\geq 2$. The class represented by $\tau$ in the homotopy group $\pi_r B(2,E)$ is non-trivial. In particular, the pull-back bundle $\tau^*(E(2,E))\rightarrow J$ is a non-trivial transitionally commutative $E$--bundle over the $r$--sphere.} 
\begin{proof}
Identifying $J$ with its nerve  we define a map of simplicial sets  
$$
\tau': J \rightarrow E(\bi,V) 
$$
by sending   an $n$--simplex $j_0\rightarrow j_2 \rightarrow \cdots \rightarrow j_n$ to the $n$--simplex
$(  j_0,  {j_1-j_0},j_2-j_1,\cdots, { j_n-j_{n-1}})$. 
Observe that there is a commutative diagram
$$
\begin{tikzcd}
J \arrow{d}{\tau'} & \sd\,J \arrow{d}{\sd\,\tau'} \arrow{l}{l} \arrow{r}{\tilde\tau} & \cC_{H(V)} \iI(V) \arrow{d} \\
E(\bi,V) & \sd\,E(\bi,V) \arrow{r}{\theta}  \arrow{l}{l} & \cC_V\iI(V)
\end{tikzcd}
$$
where $\theta$ is introduced in \ref{subtheta} and the map between the coset posets is the universal covering map.
To see that the right-hand square commutes we regard the maps as maps of posets. A chain $j_0<j_1<\cdots <j_s$ is sent to $j_0+\Span{j_1-j_0,j_2-j_0,\cdots,j_s-j_{0}}$ under the composition of $\tilde \tau$ and the covering map. On the other hand under the composition $\theta\,\sd\, \tau'$ it is sent to $j_0+\Span{j_1-j_0,j_2-j_1,\cdots,j_s-j_{s-1}}$. Note that the resulting cosets are equal.
Now $[\tau']$ is non-trivial since $[\tilde\tau]$ is non-trivial by Lemma \ref{nontrivial}. Here we have used the fact that the last vertex map is a weak equivalence, and $J$ is simply connected.
Finally, from the commutative diagram
$$
\begin{tikzcd}
J \arrow{r}{\tau'} \arrow{d}{\tau} & E(\bi,V) \arrow{d} \\
B(2,E) \arrow{r} & B(\bi,V)
\end{tikzcd}
$$
we see that $\tau:J\rightarrow B(2,E)$ represents a non-trivial class in the homotopy group.

\end{proof} 
 
Other examples of non-trivial bundles can be obtained by considering the action of the automorphism group of $E$ on the space $B(2,E)$. Equivalently, one can consider the action of $\Sp(V)$ on the coset poset $\cC_{H(V)}\iI(V)$. 
The representation of $\Sp(V)$ on the top dimensional homology of this coset poset could  be an interesting representation of the symplectic group. 
 
\subsection{Group cohomology}
We first observe a direct consequence of Theorem \ref{main2}. Then we consider the natural map $H^*(G,\FF_p)\rightarrow H^*(B(2,G),\FF_p)$ in relation to Quillen's $F$--isomorphism theorem.

\Cor{\label{cohomology} Let $R$ denote a commutative ring, and $E$   an extraspecial $p$--group. Then there is an isomorphism
$$
H^i(B(2,E),R) \cong H^i(B\pi,R)\;\; \text{for}\;\; i<r.
$$  
}
\Proof{When $r=1$ the universal cover of $B(2,E)$ is contractible since $\cC_V \iI(V)$ is a wedge of spheres (section \ref{r1}). Hence we have an isomorphism for all $i$ in this case.  If $r\geq 2$ then  the result follows from the Serre spectral sequence of the fibration in Theorem \ref{main2}.
}  

Quillen's $F$--isomorphism theorem \cite{Q71} implies that the kernel of the natural map
\begin{equation}\label{F_iso}
H^*(BG,\FF_p) \rightarrow \text{lim}_{A\in \aA(G)} H^*(BA,\FF_p)
\end{equation}
has nilpotent kernel. This map factors through the cohomology of $B(2,G)$.
 Therefore the kernel of the natural map
$$
H^*(BG,\FF_p) \rightarrow H^*(B(2,G),\FF_p)
$$
is nilpotent \cite[Proposition 3.2]{ACT12}. 

  There are two types of extraspecial $2$--groups $E^+$ and $E^-$ (see \cite[\S 23.14]{As86}).  The group $E^+$ is a central product of dihedral groups $D_8$. The cohomology of $E^+$ is detected on elementary abelian subgroups, that is, the map in \ref{F_iso} is injective. This forces   the natural map
$$
H^*(BE^+,\FF_2) \rightarrow H^*(B(2,E^+),\FF_2)
$$
to be injective as well. But this is not true in general. The quaternion group $Q_8$ is of type $E^-$, and we will show that in this case this map is not injective. 

The cohomology  ring of $Q_8$ is given by
$$
H^*(BQ_8,\FF_2)= \FF_2[x,z]/(x^2+xz+z^2,xz^2+x^2z)\otimes \FF_2[w_4]
$$
where $x,z$ are the duals of the generators of the central quotient, and $w_4$ is a $4$-dimensional class.
Next we compute the cohomology ring of $B(2,Q_8)$. In general, the cohomology ring of $B(2,Q_{2^q})$ is computed in \cite{CV17}. 
 
\Pro{\label{q8}The cohomology ring of $B(2,Q_8)$ is given by
$$
H^*(B(2,Q_8),\FF_2)= \FF_2[x,y,z]/(xy,xz,yz,x^2+y^2+z^2 )\otimes \FF_2[w_2]
$$
where $x,y,z$ are of degree one and $w_2$ is a class of degree two. As a consequence the natural map
$$
H^*(BQ_8,\FF_2)\rightarrow H^*(B(2,Q_8),\FF_2)
$$
is   not injective. 
}  
\begin{proof}
The quadratic group $Q_8$ has three maximal abelian subgroups isomorphic to $\ZZ/4$ all intersecting at the center which is a cyclic group of order two.
The space $B(2,Q_8)$ is weakly equivalent to the classifying space
$B(\ZZ/4\ast_{\ZZ/2} \ZZ/4 \ast_{\ZZ/2} \ZZ/4 )$ of the amalgamated product of the maximal abelians along the center. Consider the map between the central extensions
$$
\begin{tikzcd}
1 \arrow{r}   & \Span{w} \arrow{r} \arrow[d,equal] & \ZZ/4\ast_{\ZZ/2} \ZZ/4 \ast_{\ZZ/2} \ZZ/4  \arrow{d} \arrow{r} & \Span{x}\ast \Span{y} \ast \Span{z} \arrow{d} \arrow{r} & 1   \\
1 \arrow{r}   & \Span{w} \arrow{r}   & Q_8    \arrow{r} & \Span{x}\times \Span{z} \arrow{r} & 1  
\end{tikzcd}
$$
where $x\mapsto x$, $y\mapsto x+z$, and $z\mapsto z$.  Here $x$, $y$, $z$ are the images of the generators of the cyclic groups of order four and
each generates  a copy of $\ZZ/2$.
 In the Serre spectral sequence of the lower extension under the differential $d_2$ of the $E_2$--page,  $w$ maps to the $k$--invariant $x^2+xz+z^2$, where we denote by the same letter the dual generator in the cohomology ring. The map between the cohomology rings of the quotient groups is given by 
$$
\epsilon:\FF_2[x,z] \rightarrow \FF_2[x,y,z]/(xy,xz,yz)
$$
where $x\mapsto x+y$ and $z\mapsto z+y$. Then in the spectral sequence of the upper extension under $d_2$ the class $w$  maps to $\epsilon(x^2+xz+z^2)= x^2+y^2+z^2$ by naturality and 
$$
d_2(\Sq^1w) = \Sq^1d_2(w) = \Sq^1(x^2+y^2+z^2) =0.
$$ 
Therefore the spectral sequence collapses at the $E_3$--page and we obtain the desired result.
The map between the cohomology rings can be described explicitly as
$$
\frac{\FF_2[x,z]}{(x^2+xz+z^2,xz^2+x^2z)}\otimes \FF_2[w_4]\rightarrow \frac{\FF_2(x,y,z)}{(xy,yz,xz,x^2+y^2+z^2)}\otimes \FF_2[w_2]
$$
where the map on the first factor is induced by $\epsilon$ and $w_4$ maps to $w_2^2$.   Note that $x^2z$ maps to zero. 
\end{proof}

For odd primes $p$ it is tempting to ask how close   the map $H^*(BE,\FF_p)\rightarrow H^*(B(2,E),\FF_p)$ is to being injective? Note that  cohomology rings of extraspecial $p$--groups are not completely known for $p>2$.  
  
\subsection{Almost extraspecial $2$--group} There is another type of extraspecial group when $p=2$, also referred to as the complex type extraspecial group. It is defined to be the central product $E'=E\circ \ZZ/4$ along $\ZZ/2$. 
It sits in a central extension 
$$
0\rightarrow \ZZ/2 \rightarrow E' \rightarrow V' \rightarrow 0
$$  
where $V'=V\oplus V_0$ is isomorphic to $(\ZZ/2)^{2r+1}$. The commutator induces a bilinear form $\bi'$ on $V'$.
We have $V_0=(V')^\perp$ with respect to $\bi'$. The projection map $q:V'\rightarrow V$ induces a map of posets
$$
\hat q:\cC_{V'} \iI(V') \rightarrow \cC_V \iI(V)
$$
where the fibers $\underf{q}{v+I}$ are contractible. As a consequence of Theorem \ref{QuillenThA} this map is a weak equivalence.
Therefore the universal cover of $B(2,E')$ is a wedge of $r$--spheres. 
Analogously we can define $B(\bi',V')$ with respect to $\bi'$. There is a fibration sequence
$$
B\ZZ/2 \rightarrow B(2,E') \rightarrow B(\bi',V')
$$ 
induced by the map $E'\rightarrow V'$.   The natural map $H^*(BE',\FF_2)\rightarrow H^*(B(2,E'),\FF_2)$ can be shown to be injective using Quillen's computation in \cite{Q71a}.

\bibliography{bib}{}
\bibliographystyle{plain}

\end{document}